\documentclass[11pt]{article}
\usepackage{jmlr2e_draft}

\usepackage{amsmath, amssymb} 
\usepackage{stmaryrd}
\usepackage{epsfig,calc, graphicx}
\usepackage{wasysym}
\usepackage[usenames]{color}
\usepackage{url}
\usepackage{enumitem}
\usepackage[lined,algoruled,algosection]{algorithm2e}
\usepackage{ifthen}
\usepackage{bbm}
\usepackage{comment}
\usepackage{subcaption}
\usepackage{mathtools}
\usepackage{hyperref}
\usepackage{scalerel}
\usepackage{xcolor}

\usepackage{array}
\usepackage{booktabs}

\newcommand\ip[2]{\langle #1, #2\rangle}
 \renewcommand\sp{{n}}
\newcommand\Tin{{m}}
\newcommand\Irv{{S}}
\newcommand\Is{{I}}
\newcommand\Js{{J}}

\newcommand\indi{\mathbbm{1}}
\newcommand\one{\mathbf{1}}

\newcommand\eps{\varepsilon}

\newcommand\dico{\Phi}
\newcommand\atom{\phi}
\newcommand\pdico{\Psi}
\newcommand\patom{\psi}

\newcommand{\tp}{{t}}

\newcommand{\R}{{\mathbb{R}}}

\newcommand{\E}{{\mathbb{E}}}
\newcommand{\Imat}{{\mathbb{I}}}
\renewcommand{\P}{{\mathbb{P}}}
\newcommand{\PB}{{\P}}

\newcommand{\Dspi}{D_{\scaleto{ \sqrt{\pi}\mathstrut}{6pt}}}
\newcommand{\hol}{H}
\newcommand\DI{D_{\Is}}

\begin{document}

\title{Bounds for matrices of inclusion probabilities \\in rejective sampling}

\author{\name Simon Ruetz \email firstname.lastname@uibk.ac.at \\
	\name Karin Schnass\\
	\addr  
	University of Innsbruck\\
	Technikerstra\ss e 13\\
	6020 Innsbruck, Austria}

\editor{}

\maketitle

\begin{abstract}
Motivated by a desire to model $\sp$-sparse signals more realistically, we study matrices of inclusion probabilities for conditional Poisson or rejective sampling in the non-asymptotic regime. In particular, we provide bounds in the semi-definite ordering for matrices that collect first and second order inclusion probabilities as their diagonal and off-diagonal entries respectively, and operator norm bounds of their Hadamard products with other matrices. Finally, we demonstrate the usefulness on these bounds, which only involve first order inclusion probabilities, in a toy application from dictionary learning.
\end{abstract}



\section{Introduction}
The goal of this paper is to provide bounds for matrices of inclusion probabilities in rejective sampling and thus make rejective sampling a practicable model for $\sp$-sparse signals. Rejective sampling, known also as conditional Poisson/Bernoulli sampling, selects a subset $\Irv$ of size $\sp$ from $[K] =\{1,\ldots, K\}$, by flipping $K$ independent Bernoulli switches $\delta_k$ with $\P(\delta_k = 1) = p_k \in (0,1)$ and keeping the set of active indices $\Irv=\{k:\delta_k=1\}$ if it has the correct size. First order inclusion probabilities then count how often such a set contains a fixed index $i$, $\pi_i(\sp) = \P(\delta_i =1 | \sum_k \delta_k = n) = \P_n(i\in \Irv)$, and higher order inclusion probabilities how often they contain a fixed set $L = \{i_1, \ldots, i_\ell\}$, $\pi_L(\sp) = \pi_{i_1\ldots i_\ell}(n) = \P_n(L\subseteq \Irv)$. The bounds presented here are of the following form.
\begin{ftheorem}
Consider rejectively sampling sets of size $\sp$ with generating probabilities~$p$. Let $\pi=\pi(\sp)$ be the vector of first order inclusion probabilities, $D_{\pi}$ the diagonal matrix with $\pi$ on the diagonal and $\Sigma_\sp$ be the matrix of first and second order inclusion probabilities, ie. with $ij$-th entry $\pi_{ij}(\sp)$ and $\pi(\sp)$ on the diagonal. Further, let $\preceq$ denote the semi-definite order and $\odot$ the Hadamard (entrywise) product.
Then we have 
\begin{align*}
 \Sigma_\sp  &\preceq \pi \pi^\tp + 2 D_{\pi} \quad \mbox{as well as}\\
	\| A \odot \Sigma_\sp \| & \leq \kappa \cdot \|D_{\pi} [A -A\odot \Imat \,] D_{\pi}\| + \| A \odot D_{\pi}\|,  
\end{align*} 
for $A\in \R^{K\times K}$ and $\kappa = (1+\|\pi\|_\infty) (1-\|\pi\|_\infty)^{-2}$. 
\end{ftheorem}
Since the bounds above are valid in the non-asymptotic setting, they make rejective sampling usable and useful for people working in sparse signal processing, data science or machine learning, where population size~$K$ and sample size~$\sp$ cannot be assumed to be large. Indeed, our main message is, that if you have an equality in uniform $\sp$-sampling or Poisson sampling that features $\sp/K$ or $p$, you can often get an analogue inequality for rejective sampling, featuring instead $\pi$ or even $p$ with a benign constant.\\
To motivate the renewed interest in a classical topic such as inclusion probabilities of rejective sampling, we next provide a motivation for and introduction to rejective sampling aimed at practitioners and interested sampling experts using an application example in the modelling of sparse signals. Section~\ref{sec_main_results} contains the main result in more detail together with some useful auxiliary results. We also provide a short application example in dictionary learning, while a more elaborate example can be found in Appendix~\ref{sec_dico_update}. Finally, Section~\ref{sec_proofs} contains the proofs of the main result. 

\subsection{Motivation and rejective sampling model}
One of the most useful priors in signal processing, e.g. for signal acquisition and restoration such as compressed sensing, denoising or inpainting, \cite{do06cs, donoho:bp, maelsa08}, is that the signals of interest are sparse in a dictionary or basis. This means that there exists a set of vectors $\atom_i\in \R^d$, called atoms and collected in the dictionary $\dico = (\atom_1, \ldots, \atom_K)\in \R^{d\times K}$, such that we can write each signal $y$ of interest as
\begin{align*} 
y = \dico x = \sum_{k=1}^K  x(k) \atom_k = \sum_{k\in \Is } x(k) \atom_k,
\end{align*}
where the coefficients $x\in \R^K$ are $\sp$-sparse, which means that only $\sp$ entries $x(k)$ are non-zero. The (sorted) set $\Is=\{k : x(k) \neq 0\}$, which collects all indices $k$ of non-zeros entries $x(k)$, is then called the support of the sparse signal and sometimes, it will be convenient to rewrite the sparse signals as superposition of scaled dictionary atoms in the support in several ways, meaning for $\dico_\Is = (\atom_{i_1}, \ldots, \atom_{i_\sp})$ and $x_\Is=(x_{i_1},\ldots, x_{i_\sp})$, where we usually assume that $i_j \leq i_{j+1}$, we have
\begin{align*} 
y  = \dico_\Is x_\Is = \dico D_x \one_\Is = \dico D_{\Is} x .
\end{align*}
Here $\one_\Is$ denotes the indicator vector of the set $\Is \subseteq [K]=\{1,\ldots ,K\}$, meaning $\one_\Is(j) = 1$ for $j\in\Is$ and zero else, while $D_x$ for $x \in \R^K$ and $D_{\Is}$ for $\Is \subseteq [K]$ denote the diagonal matrices with $x$ and respectively $\one_\Is$ on the diagonal.
\\
Ideally, given an algorithm and the sparse model, meaning the dictionary $\dico$ and the sparsity level $\sp$, one can quantify whether the algorithm is successful. For instance if the dictionary has small coherence $\mu(\dico) = \max_{i\neq j} |\ip{\atom_i}{\atom_j}|$, $\ell_1$-minimisation based schemes, \cite{chdosa98,donoho:bp}, or OMP, \cite{parekr93}, can recover the sparse support and coefficients $x_\Is$ from $y$ as long as $2\mu \sp = 2\mu |\Is| <1$, meaning one can compress $y$ to $x_{\Is}$, \cite{bp:fuchs,Tropp:greed}. \\
Similarly, if one has the design freedom to choose $\dico$, the task of reconstructing any $\sp$-sparse $x\in \R^K$ from $y=\dico x \in \R^d$ is known as compressed sensing. There one of the key results is that if one chooses $\dico$ randomly, e.g. with i.i.d Bernoulli entries, then as long as $ \sp \cdot \log(K) / d \lesssim 1$ with high probability all $\sp$-sparse $x\in \R^K$ can be reconstructed via $\ell_1$-minimisation, \cite{cata06, carota06}. Since for $K > d$ we have $\mu$ at best of order $1/\sqrt{d}$, choosing a random dictionary, means that the scheme works for much larger $\sp$. \\
However, in many situations the dictionary $\dico$ is predetermined by the application, meaning it cannot be chosen or modelled as random. At the same time all deterministic coherence based results are too weak to be useful or do not reflect the practical performance of an algorithm well. \\ 
In these situations one can try to get a better fit between theory and practical performance by adopting a random model on the signal coefficients, meaning the sparse support and the coefficient size, and by studying the average performance. After all, most algorithms are designed to be used for many signals or in the case of machine learning they use as input a batch of \emph{typical} signals.\\
Going back to our example of sparse approximation, one can show that if the signals are modelled as
\begin{align} \label{random_sparse_model}
y = \sum_{k\in \Irv } c(k) \sigma(k) \atom_k = \dico_\Irv (c_\Irv \odot \sigma_\Irv),
\end{align}
where $\Irv$ chooses supports $\Is$ uniformly at random among all sets of size $\sp$, the signs~$\sigma$ are a Rademacher sequence, $c$ is a strictly positive sequence and $\odot$ is the Hadamard (entrywise) product, then with high probability one can recover the support and non-zero coefficients via $\ell_1$-minimisation, as long as $\sp \mu^2 \log(K) \lesssim 1$, \cite{tr08}.\\
In order to reflect that some atoms usually appear with higher coefficients than others, it is of course easy to add a random model on $c$, e.g. by letting each entry $c(k)$ follow a distribution bounded away from zero (and bounded) but with different mean.
\begin{figure}[t] 
  \centering
  \includegraphics[width=\linewidth]{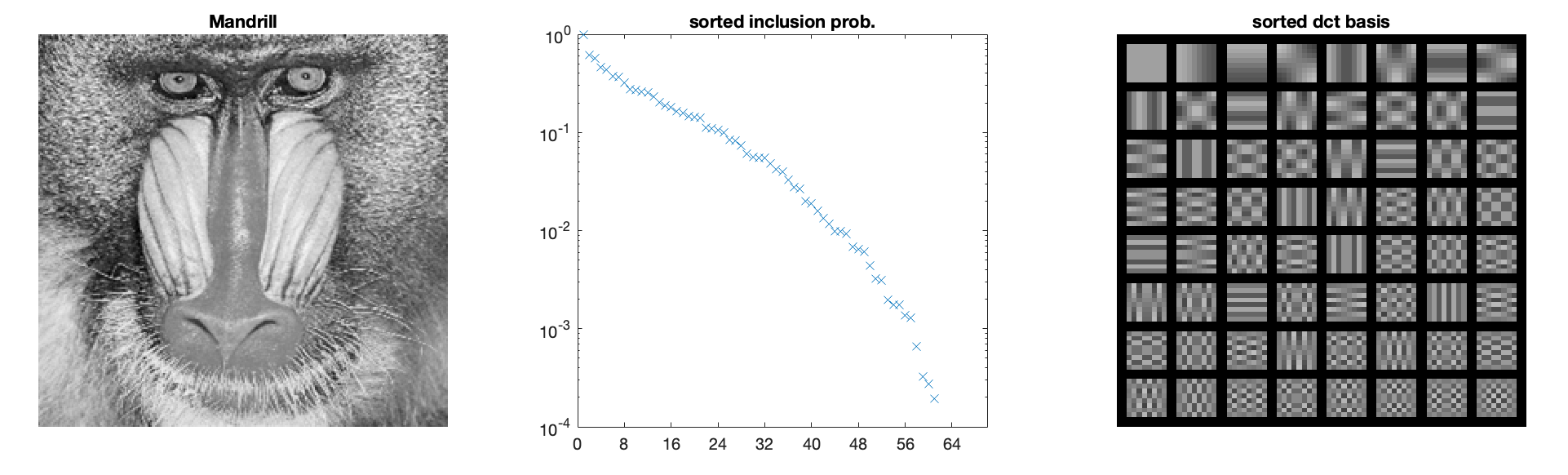}
  \caption{For each $8\times8$-patch in \emph{Mandrill} (left) we find the best approximation using $8$ atoms of the DCT basis (right). We then count how often each atom is used and calculate the empirical inclusion probabilities (middle). The DCT atoms are sorted from most to least frequent row-wise, i.e. most frequent top left, least frequent bottom right.  \label{mandrill_incl_prob}}
\end{figure}
However, it is a little more tricky to incorporate another property that we are likely to encounter in a realistic set of sparse signals, i.e. that some atoms are more likely to be in the $\sp$-sparse support that others.  For instance, if we look at the best approximation of all $8\times 8$ patches of \emph{Mandrill} with $\sp= 8$ DCT atoms, and count how often each DCT atom is used, we see in Figure~\ref{mandrill_incl_prob} that some atoms clearly appear for frequently than others. \\
In order to get a still simple but more realistic model, let us forget for a moment that we want only $\sp$ atoms in the support. In this case we can give each atom a probability to be in the support, meaning we define a series of independent Bernoulli 0-1 variables $\delta_k$ which are 1 with probability $p_k \in (0,1) $ and zero with probability $(1-p_k)$. 
Then for the random vector $\delta$ and each $v \in \{0,1\}^K=\Delta$ we have 
$$\PB(\delta = v) = \prod_{i:v(i) =1 }p_i\prod_{j:v(j) = 0}(1-p_j)=: p^v \cdot (\one - p)^{1-v}.$$ 
Now if we identify vectors in $\Delta$ with subsets of $[K]$ via $\iota : \Delta \to \mathcal{P}([K]), v\mapsto \{ i : v(i) =1\}$, meaning $\iota^{-1}(\Is) = \one_{\Is}$ we get the set-valued random variable $\Irv = \iota \circ \delta$, which is distributed as 
\begin{equation*}
    \PB(\Irv = \Is) := \prod_{i \in \Is}p_i\prod_{j \notin \Is}(1-p_j) =: p^\Is (\one - p)^{\Is^c} \quad \mbox{for}  \quad \Is \subseteq [K].
\end{equation*}
Choosing a support or sample according to the distribution above is called Poisson sampling. 
Now if our probability weights $p_k$ sum to $\sp$, meaning for the vector $p = (p_1, \ldots ,p_K)$ we have $\|p\|_1 =\sp$,
we know that in expectation the size of the chosen supports is $\sp$, since
\begin{align*}
	\E(| \Irv | ) = \E \Big(\sum_{k=1}^K \one_\Irv(k)\Big) =  \E \Big(\sum_{k=1}^K\delta_k \Big) = \sum_{k=1}^K p_k = \|p\|_1 = \sp,
\end{align*}
and concentration of measure tells us that most supports that we choose with this method will have support size close to $\sp$. Unfortunately, this does not mean that most supports have actual size $\sp$. However, we we can force this property by conditioning on it, meaning we define a new distribution of $S$ via
\begin{equation*}
    \P_\sp(\Irv = \Is) := \PB(  \Irv = \Is \; \mid \;  |\Irv| = \sp) = \begin{cases}
               c^{-1}\cdot  p^\Is (\one - p)^{\Is^c} \quad& \mbox{if} \quad |\Is|=\sp\\
               0 \quad& \mbox{else}
            \end{cases},
\end{equation*}
where $c = \PB(|\Irv|=\sp)= \sum_{\Is:|\Is| = \sp}  p^\Is (\one - p)^{\Is^c} $.
Due to its construction, choosing a support or sample according to the distribution above is called conditional Poisson sampling, conditional Bernoulli sampling or rejective sampling, since we reject all samples of the wrong size. Conditional Poisson sampling was popularised and intensively studied by J.~H\'ajek, e.g. \cite{hajek1964, hajek1981}. It belongs to a larger class of sampling measures called exponential designs and we refer to~\cite{tille06} for a thorough introduction, \cite{Aires99} for algorithmic aspects and \cite{Yam12, Ber16} for more recent developments. \\
Coming back to our goal of modelling $\sp$-sparse signals, we can now generalise the model in~\eqref{random_sparse_model}, by allowing the $\sp$-sparse supports to be chosen according to rejective sampling with parameter vector $p$, such that $\|p\|_1 =\sp$. If we then consider our sparse approximation problem, indeed we get an analog statement, that with high probability one can recover the sparse supports via $\ell_1$-minimisation, as long as the hollow Gram matrix of the dictionary, $H=\dico^\tp \dico-\Imat_K$, weighted by the diagonal matrix $\Lambda$ with $\Lambda_{ii} =\sqrt{p_i}$, has small operator and column-wise Euklidean norm, that is $\|\Lambda H \Lambda\|_{2,2}\lesssim 1$ and $\|H\Lambda\|_{\infty,2}^2 \log(K) \lesssim 1$, \cite{rusc21}.\\
While rejective sampling clearly allows to model $\sp$-sparse signals more accurately, it comes with certain complications or drawbacks. For instance, the
first order inclusion probabilities $\pi_k= \pi_k(\sp)= \P_\sp(k \in \Irv)$, that we can observe, differ from the generating Bernoulli probabilities $p_k=\PB(k \in \Irv)$, we use in the model, meaning unless $p_k = c$ for all $k$, we have 
\begin{align*}
\pi = \pi(\sp) = \E_\sp[\one_{\Irv}] \neq \E[\one_{\Irv}] = p.
 \end{align*}
Similarly, the higher order inclusion probabilities are not a product of the first order inclusion probabilities, so for $L\subseteq [K]$ we have
\begin{align*}
\pi_L(\sp)= \P_\sp(L \subseteq \Irv) \neq \prod_{\ell \in L} \pi_\ell \quad \mbox{vs.} \quad \prod_{\ell \in L} p_\ell = \PB(L \subseteq \Irv).
\end{align*}
This also means that the matrix of inclusion probabilities $\Sigma_\sp = \E_\sp [\one_{\Irv} \one_{\Irv}^\tp]$ does not have a simple structure
like $\Sigma =  \E [\one_{\Irv}\one_{\Irv}^\tp] = pp^\tp + D_{p}(\Imat-D_p)$, so weighting a matrix $A$ with $\Sigma_\sp$, is not as straightforward as weighting with $\Sigma$, where we simply have
\begin{align*}
A\odot \Sigma = D_p A D_p +  A \odot D_p(\Imat-D_p).
\end{align*}
However, the results presented in the next section show that many equalities in uniform $\sp$-sampling or Poisson sampling that feature $\sp/K$ or $p$ respectively, have analogue inequalities in rejective sampling, featuring instead of $\sp/K$ or $p$ either $\pi$ or even $p$ with a benign constant. As such they provide the key to making rejective sampling an attractive and usable sparse model.

\section{Main results}\label{sec_main_results}
We start with a simple comparison bound between first order inclusion probabilities $\pi$ and generating probabilities $p$, which will later allow us to transfer all matrix bounds in terms of $\pi$ into bounds in terms of $p$.
\begin{lemma}[Comparison bound]\label{comparison_bound}
Let $\pi_i(\sp)= \P_\sp( i \in \Irv)$ be the first order inclusion probabilities associated to a rejective sampling model with parameter $\sp$ and weight vector $p=(p_1, \ldots, p_K)$ with $p_i \in (0,1)$ and $\|p\|_1=\sp$, then we have
\begin{align*}
1-\|p\|_\infty  &\leq \frac{\pi_i(\sp)}{p_i}  \leq 2 . 
\end{align*}
 \end{lemma}
\noindent The relation between $\pi_i(\sp)$ and $p_i$ was first studied by H\'ajek, \cite{hajek1964, hajek1981}, who showed that
\begin{align*}
    \max_{1 \leq i \leq K} \left| \pi_i / p_i - 1 \right| \rightarrow 0 \quad \text{as} \quad q:= \sum_{i = 1}^K p_i (1 - p_i)  \rightarrow \infty,
\end{align*}
and conjectured, as later shown by Milbrodt, \cite{milbrodt92}, that for $\alpha_i := \pi_i(1) = \frac{p_i}{1-p_i}$, 
\begin{align*}
  \min_i \alpha_i =  \min_i \frac{p_i}{1-p_i} &\leq \min_i \pi_i \qquad \text{and} \qquad \max_i \pi_i \leq \max_i \frac{p_i}{1-p_i} = \max_i \alpha_i .
\end{align*}
The result above, which already gives control over the extreme points of the two sequences, was further generalised in \cite{kochar01,Yam12} to
\begin{align}
    \pi \prec \alpha \quad \text{and} \quad 
    (K^{-1},\dots , K^{-1})= \frac{\pi(K)}{K} \prec \dots \prec \frac{\pi(\sp)}{\sp} \prec \dots \prec \pi(1)= \alpha,\label{incs1}
\end{align}
where $u \prec v$ for two vectors $u,v \in\R^K$ means, that after sorting them in increasing order, their sorted versions $\vec{u},\vec{v}$ satisfy
\begin{align*}
    \sum_{k= j}^K \vec{u}_k \leq \sum_{k = j }^K \vec{v}_{k}, \qquad \forall j =1, 2,\dots, K.
\end{align*}
The advantage of the bound here is that, while maybe less sharp, it holds in the non-asymptotic setting, where $q\leq S \ll K$, and that it is valid for all pairs $(p_i,\pi_i)$, rather than just the extremal ones or partial (sorted) sums.\\
While the comparison bound above has its asymptotic analogue, bounds in the operator norm or semi-definite order for matrices of inclusion probabilities, such as $\Sigma_\sp$ with $\pi_i(n)$ and $\pi_{ij}(n)$ as diagonal and off-diagonal entries, seem to be a unstudied so far. However, $\Sigma_\sp$ appears quite naturally in a sparse context. \\
For instance, the idea behind many dictionary learning algorithms is to construct glimpses of the dictionary $\dico$ using a current guess $\pdico$. Assuming a sparse signal model as in \eqref{random_sparse_model} with constant coefficient sequence $c$ these glimpses can (in first order approximation) be of the form
$\dico D_\Irv + \dico D_\Irv  Z^\tp \pdico D_\Irv $ for $Z = \dico -\pdico$, so that aggregating them one arrives at
\begin{align}
\hat \pdico:= \dico\, \E[  D_\Irv  + D_\Irv  Z^\tp \pdico D_\Irv  ] = \dico [D_\pi +  (Z^\tp \pdico) \odot \Sigma_\sp]. \label{eq_newdico}
\end{align}
Determining whether a normalised version of $\hat \pdico$ is a better approximation to $\dico$ than $\pdico$, then becomes a question of bounding 
$(Z^\tp \pdico) \odot \Sigma_\sp$ or rather bounding $H \odot \Sigma_\sp$ where $H =  Z^\tp \pdico - ( Z^\tp \pdico)\odot \Imat$.
%
\begin{theorem}[Operator norm \& semi-definite order bounds] \label{matrix_bounds}
Let $\E_\sp$ be the expectation according to the rejective sampling probability with parameter $\sp$ and weights $p_i \in (0,1)$. Further let 
$\pi\in \R^K$ be the vector of first order inclusion probabilities, meaning $\pi = \E_{\sp}[\one_\Irv]$, $D_\pi$ be the $K\times K$ matrix with $\pi$ on the diagonal and zero else, and $\Sigma_\sp$ be the matrix of (up to) second order inclusion probabilities, meaning $\Sigma_\sp =\E_\sp[\one_\Irv  \one_\Irv^\tp]$. Then for the Hadamard (entrywise) product with any matrix $A\in \R^{K\times K}$ we have 
\begin{align} \label{matrix_a}
   \| A \odot \Sigma_\sp\|   \leq \frac{1 + \|\pi\|_\infty}{(1 - \|\pi\|_\infty)^2} \cdot \|D_{\pi} [A -A\odot \Imat \,] D_{\pi}\| + \| A \odot D_{\pi}\|. \tag{a} 
\end{align}

\noindent Writing $A\preceq B$ for two positive semi-definite (p.s.d) matrices if $A-B$ is p.s.d, we have
\begin{align}\label{matrix_b}
\Sigma_\sp  &\preceq \pi \pi^\tp + 2 D_{\pi}, \tag{b}
\end{align} 
and, defining for $L\subseteq [K]$ with $|L|<\sp$ the set $\mathcal L: = \{\Is\subseteq [K] : L\subseteq \Is\}$, we have
 \begin{align} \label{matrix_c}
\E_{\sp} \big[ \one_{\Irv\setminus L} \one_{\Irv\setminus L }^\tp \cdot \indi_\mathcal{L}(\Irv) \big] 
&\preceq  \Sigma_{\sp-|L|} \cdot \prod_{\ell \in L} \frac{\pi_\ell}{1-\pi_\ell}. \tag{c}
\end{align}
\end{theorem}
%
Interestingly the theorem above allows us to bound quantities involving second order inclusion probabilities using only first order inclusion probabilities or in combination with Lemma~\ref{comparison_bound} even generating probabilities, which are much easier to handle. For instance, we can answer the question how close $\hat \pdico$ is to a scaled version of $\dico$.
\begin{example}[Toy application in dictionary learning] 
Recall that we have a dictionary $\dico$, a guess $\pdico= \dico + Z$ and are wondering when the updated guess $\hat\pdico= \dico [D_\pi +  (Z^\tp \pdico) \odot \Sigma_\sp]$ will be an improvement over $\pdico$.\\
We set $\eps_i = \ip{z_i}{\patom_i}$ so that $D_\eps = ( Z^\tp \pdico)\odot \Imat$ and $H = Z^\tp \pdico - D_\eps$. We can then rewrite
$\hat \pdico = \dico [ D_\pi (\Imat + D_\eps) + H \odot \Sigma_\sp]$
and get abbreviating  $C_\pi = (1 + \|\pi\|_\infty)(1 - \|\pi\|_\infty)^{-2}$ and $ \Dspi = D_\pi^{1/2}$ that
\begin{align*}
\| \hat \pdico (\Imat + D_\eps)^{-1}\Dspi^{-1} - &\dico  \Dspi  \| =   \| \dico \Dspi \, \Dspi^{-1} \, ( H  \odot \Sigma_\sp)\, \Dspi^{-1}(\Imat + D_\eps)^{-1}\| \\
&\leq \| \dico\Dspi\| \cdot \| ( \Dspi^{-1} H \Dspi^{-1} ) \odot \Sigma_\sp \|  \cdot \|(\Imat + D_\eps)^{-1}\|\\
&\leq \| \dico \Dspi\| \cdot C_\pi \cdot \|  \Dspi H \Dspi  \| \cdot \|(\Imat + D_\eps)^{-1}\|\\
&\leq C_\pi \cdot \| \dico \Dspi\| \cdot \|(\Imat + D_\eps)^{-1}\| \cdot ( \|  Z \Dspi \| \cdot \|\pdico \Dspi \| + \|D_\eps D_\pi \| ) \\
&\leq  C_\pi \cdot \| \dico \Dspi\| \cdot \|(\Imat + D_\eps)^{-1}\| \cdot ( \|\pdico \Dspi \| + \|\Dspi\| ) \cdot \|  Z \Dspi\| .
\end{align*}
Since $\| Z \Dspi\| = \|(\dico-\pdico) \Dspi\|$ is the weighted distance between the current estimate $\pdico$ and the desired dictionary $\dico$, it is now easy to infer from the bound above under which conditions on $\dico$ and $\pdico$ the normalised version of $\hat \pdico$ or equivalently of $\hat \pdico (\Imat + D_\eps)^{-1}\Dspi^{-1}$ will be a better approximation to $\dico$ in this weighted distance. 
\end{example}
Moreover, using Theorem~\ref{matrix_bounds}, it is also possible to estimate terms involving higher order inclusion probabilities of the form $\E_{\sp} [ D_{\Irv} H D_{\Irv} H^\tp D_{\Irv}]$, where $H$ is again a matrix with zero diagonal. In dictionary learning such terms appear considering glimpses in second order approximation, and we can bound them using the following corollary of Theorem~\ref{matrix_bounds}. 
\begin{corollary}\label{cor_bound_HHt}
Under the same conditions as in Theorem~\ref{matrix_bounds} we have for any matrix $H$ with zero diagonal that
\begin{align*}
\big \| \E_{\sp} [ D_{\Irv} H D_{\Irv} H^\tp D_{\Irv} ]\big \| \leq
(1 + \|\pi\|_\infty)\cdot  \frac{\|D_\pi H D_\pi H^\tp D_\pi \|}{(1 - \|\pi\|_\infty)^3}   +  \frac{\|D_\pi \odot (H D_\pi H^\tp) \| }{1 - \|\pi\|_\infty} 
. 
\end{align*}
\end{corollary}
Note that $D_\pi \odot (H D_\pi H^\tp)$ is a diagonal matrix, so we can also write $\|D_\pi \odot (H D_\pi H^\tp) \| = \max_k \pi_k (H D_\pi H^\tp)_{k,k}$.
The interesting fact about the bound above is, that it allows to bound the matrix $\E_{\sp}[ D_{\Irv} H D_{\Irv} H^\tp D_{\Irv} ]$, which involves up to of third order inclusion probabilities, again using only the first order inclusion probabilities $\pi$. \\
For the interested reader we provide a short explanation how the first and second order glimpses in dictionary learning are constructed and how they can be bounded using Corollary~\ref{cor_bound_HHt} in Appendix~\ref{sec_dico_update}. For even more details and a real application we refer to \cite{rusc26}.\\
\noindent Here we next we turn to the proofs of both results.

\section{Proof of main results}\label{sec_proofs}
\noindent We start with the proof of Lemma~\ref{comparison_bound}.
\begin{proof}[of Lemma~\ref{comparison_bound}]
\noindent We want to show that $1-\|p\|_\infty  \leq \pi_i(\sp)/p_i  \leq 2 $ whenever $\|p\|=\sp$.
The upper bound follows directly from \cite[Lemma~7]{rusc21} so we only need to show the lower bound. We abbreviate $\pi_i = \pi_i(\sp)$ and $c: = 1-\|p\|_\infty \leq  \pi_i/p_i$ as well as $\PB(\Is) = \PB(\Irv =\Is)$ for any set $I\subseteq [K]$. By definition, we have
\begin{align*}
  \pi_i &= \PB( i \in \Irv \mid |\Irv| = \sp ) = \frac{\PB(\{ i \in \Irv\} \cap \{|\Irv| = \sp \})}{\PB(|\Irv| = \sp )} 
    = \frac{\sum_{\Is:|\Is| = \sp, i \in \Is} \PB(\Is)}{\sum_{\Is:|\Is| = \sp}\PB(\Is)} \underbrace{\sum_{\Js}\PB(\Js)}_\text{\clap{= 1}},
\end{align*}
and $p_i = \sum_{\Js : i \in \Js} \PB(\Js)$. So the desired inequality $c\cdot p_i \leq \pi_i$ is equivalent to 
\begin{align*}
 c  \sum_{\Is:|\Is| = \sp}\PB(\Is)  \sum_{\Js: i \in \Js} \PB(\Js)  \leq  \sum_{\Is:|\Is| = \sp, i \in \Is} \PB(\Is) \sum_{\Js}\PB(\Js).
\end{align*}
Splitting the sum over $\Is,\Js$ into sums over those containing $i$ and those not containing $i$ we see that the inequality above is implied by
\begin{align}
  c \sum_{\Is:|\Is| = \sp, i\notin \Is}\PB(\Is)  \sum_{\Js: i \in \Js} \PB(\Js) \leq \sum_{\Is:|\Is| = \sp, i \in \Is} \PB(\Is) \sum_{\Js: i\notin \Js}\PB(\Js). \label{p_leq_2pi_a}
\end{align}
Note that for any set $\Is$ not containing the index $i$ we have 
$$\frac{p_i}{1-p_i} \cdot \PB(\Is) = \frac{p_i}{1-p_i} \prod_{k \in \Is} p_k \prod_{k\notin \Is} (1-p_k) =  \prod_{k \in \Is\cup\{i\}} p_k \prod_{k\notin \Is\cup\{i\}} (1-p_k) =  \PB(\Is\cup\{i\}).$$
Multiplying both sides in \eqref{p_leq_2pi_a} with $p_i/(1-p_i)$ we get
\begin{align*}
  c \sum_{\Is:|\Is| = \sp+1, i\in \Is}\PB(\Is)  \sum_{\Js: i \in \Js} \PB(\Js) \leq \sum_{\Is:|\Is| = \sp, i \in \Is} \PB(\Is) \sum_{\Js: i\in \Js}\PB(\Js),
\end{align*}
so it suffices to show that
\begin{align*}
 c  \sum_{\Is:|\Is| = \sp+1, i\in \Is}\PB(\Is)  \leq \sum_{\Is:|\Is| = \sp, i \in \Is} \PB(\Is).
\end{align*}
Indeed, we have 
\begin{align*}
      c \sum_{\Is:|\Is| = \sp+1, i\in \Is}\PB(\Is)  &= c  \sum_{\Is:|\Is| = \sp+1, i\in \Is}\PB(\Is)  \sum_{k:k \in \Is, k\neq i}\frac{1}{\sp} \\
     &=c \sum_{\Is:|\Is| = \sp+1, i\in \Is} \frac{1}{\sp} \sum_{k:k \in \Is, k\neq i} \PB(\Is\setminus\{k\}) \frac{p_k}{1-p_k} \nonumber \\
     &\leq \frac{c}{\sp(1-\|p\|_\infty)} \sum_{\substack{(\Is,k):|\Is| = \sp +1,i\in \Is\\\phantom{la} k\in \Is, k\neq i}} \:  \PB(\Is\setminus\{k\}) \cdot p_k \\
     &= \frac{1}{\sp} \sum_{\Js: |\Js| = \sp, i \in \Js} \PB(\Js) \sum_{k \notin \Js}  p_k \leq \sum_{\Js: |\Js| = \sp, i \in \Js} \PB(\Js),
\end{align*}
where we used that $\sum_{k \notin \Js}  p_k \leq \sum_k p_k = \sp$. 
\end{proof}
In order to prove Theorem~\ref{matrix_bounds} we need two auxiliary results. The first is a corollary of Schur's product theorem, see e.g. \cite{horn13}. For convenience we include its short proof. 
\begin{corollary}[of Schur's theorem]\label{had} 
Let $A,B\in\R^{K\times K}$ be two square matrices of the same dimension. If $A$ is positive-semidefinite (p.s.d.), then
\begin{align*}
    \| A \odot B \| \leq \|A \odot \Imat\| \cdot \|B\| .
\end{align*}
\end{corollary}
\begin{proof}
The matrix
\begin{align*}
    \begin{pmatrix}
     \|B\| \cdot (A \odot \Imat) & A \odot B \\
     (A \odot B)^\tp & \|B\|\cdot (A \odot \Imat)
    \end{pmatrix} = 
    \begin{pmatrix}
     A & A \\
     A  & A
    \end{pmatrix} \odot 
    \begin{pmatrix}
     \|B\| \cdot \Imat & B \\
     B^\tp & \|B\| \cdot \Imat
    \end{pmatrix}
\end{align*}
is p.s.d., since the right hand side of the equation is a Hadamard product of two p.s.d. matrices and thus by Schur's product theorem also p.s.d. Further by~\cite[Theorem 7.7.9]{horn13} this means that there exists a contraction $C$, meaning $\|C\| \leq 1$, such that
\begin{equation*}
    A \odot B = \|B\| \cdot (A \odot \Imat)^{1/2} \: C \: (A \odot \Imat)^{1/2},
\end{equation*}
and hence $\|A \odot B\| \leq \|A \odot \Imat\| \cdot \|B\| $.
\end{proof}
The second result is the following lemma, which collects relations between first and second order inclusion probabilities for $\sp$ resp. $\sp -1$, corresponding to \cite[Eq. 5.21]{tille06}, and bounds between first order inclusion probabilities for $\sp$ resp. $\sp -1$. 
\begin{lemma}\label{scalar_bounds}
Let $\pi_i(\sp)= \P_\sp( i \in \Is)$ be the first order inclusion probabilities associated to a rejective sampling model with parameter $\sp$ and weight vector $p=(p_1, \ldots, p_K)$ with $p_i \in (0,1)$, then we have
\begin{align}
 & \pi_{ij}(\sp) =  \pi_i(\sp) \cdot \frac{\pi_j(\sp-1) - \pi_{ij}(\sp-1)}{1 - \pi_i(\sp-1)}  \quad \mbox{for all } i\neq j  \label{tille_521}\\
& \mbox{and}\quad   \pi_i(\sp-1) \leq  \pi_i(\sp) \quad \mbox{for all } i.  \label{scalar_b}
\end{align}
\end{lemma}
A proof of Eq.~\eqref{tille_521} can be found in \cite{tille06}. The inequality in \eqref{scalar_b} follows from \cite{neg_dep_09}, using that $\mu=\PB$  is strongly Raleigh. Therefore $\mu_{n-1} = \P_{n-1} \preceq \P_n =\mu_n$, which further implies that for increasing functions $f$,  such as $f(S) = \indi_{\Irv}(k)$, $\pi_i(\sp-1) = \int f \mathrm{d}\mu_{\sp-1} \leq  \int f \mathrm{d}\mu_{\sp} =\pi_i(\sp)$. For convenience, alternative proofs, which use only the language and techniques introduced in this paper, are provided in Appendix~\ref{sec_proof_scalar_bounds}.
With the two results above we are ready to prove Theorem~\ref{matrix_bounds}.
\begin{proof}[of Theorem~\ref{matrix_bounds}]
\noindent \textbf{(a)} We want to show that 
\begin{align*}
\| A \odot \Sigma_{\sp} \|  \leq \frac{1 + \|\pi\|_\infty}{(1 - \|\pi\|_\infty)^2} \cdot \|D_{\pi} [A -A \odot \Imat\,] D_{\pi}\| + \| D_{\pi} \odot A \|.
\end{align*}
First note that since $\Sigma_{\sp}$ has $\pi(\sp)$ on the diagonal, a simple application of the triangle inequality yields
\begin{align*}
    \| A \odot \Sigma_{\sp} \| &\leq \| A \odot \Sigma_{\sp} - A\odot D_{\pi(\sp)} \|  +\|  A \odot D_{\pi(\sp)}\|=  \| (A -A \odot \Imat) \odot \Sigma_{\sp} \| + \| A\odot D_{\pi(\sp)} \|,
\end{align*}
which already proves the theorem for $\sp=1$, where all off-diagonal entries of $\Sigma_{\sp}$ are zero, meaning the first norm term is zero.
In case $\sp\geq 2$ it remains to show that for $\hol= A-A \odot \Imat$ we have $ \| \hol \odot \Sigma_{\sp} \| \leq c \cdot \| D_{\pi(\sp)} \hol  D_{\pi(\sp)} \|$ with constant $c$ as above. Using the abbreviation $\Delta_{\pi(\sp)} = \Imat - D_{\pi(\sp)}$ we know from \eqref{tille_521} that 
 \begin{align*}
(\hol \odot \Sigma_{\sp} )_{ij} &= \hol_{ij} \cdot \pi_{ij}(\sp)= \frac{ \pi_i(\sp) }{ 1 - \pi_i(\sp-1)} \cdot H_{ij} \cdot \left[ \pi_j(\sp-1) - \pi_{ij}(\sp-1)\right]\\
&= \left( \Delta_{\pi(\sp-1)}^{-1} D_{\pi(\sp)}  \hol  D_{\pi(\sp-1)}- \Delta_{\pi(\sp-1)}^{-1} D_{\pi(\sp)} \hol \odot \Sigma_{\sp-1} \right)_{ij}.
\end{align*}
Next Lemma~\ref{scalar_bounds} tells us that $D_{\pi(\sp-1)}\preceq D_{\pi(\sp)}$, which leads to
 \begin{align}
\| \hol \odot \Sigma_{\sp} \| &\leq \| \Delta_{\pi(\sp-1)}^{-1}\|\cdot\| D_{\pi(\sp)} \hol  D_{\pi(\sp-1)}-  D_{\pi(\sp)}  \hol \odot \Sigma_{\sp-1}  \| \notag \\
 &\leq (1-\|\pi(\sp)\|_\infty)^{-1}\cdot \left( \|D_{\pi(\sp)}  \hol  D_{\pi(\sp)}\| + \|D_{\pi(\sp)} \hol \odot \Sigma_{\sp-1} \| \right). \label{bound_D_pi_left}
 \end{align}
Applying the inequality above to $H^\tp$ and using the symmetry of $\Sigma_{\sp}$ we also get
\begin{align}
 \| \hol \odot \Sigma_{\sp} \| &\leq (1-\|\pi(\sp)\|_\infty)^{-1}\cdot \left( \|D_{\pi(\sp)}  \hol  D_{\pi(\sp)}\| + \| \hol D_{\pi(\sp)} \odot \Sigma_{\sp-1} \| \right). \label{bound_D_pi_right}
 \end{align}
For $\sp=2$, the matrix $\Sigma_{\sp-1} $ is again a diagonal matrix, meaning the second norm term vanishes and we are done. For $\sp>2$ we simply apply the inequality in \eqref{bound_D_pi_right} to $\bar \hol \odot \Sigma_{\sp-1} $ with $\bar \hol = D_{\pi(\sp)} \hol$, leading to 
\begin{align*}
 \| D_{\pi(\sp)} \hol &\odot \Sigma_{\sp-1} \| = \|  \bar \hol \odot \Sigma_{\sp-1} \|\\
 &\leq (1-\|\pi(\sp-1)\|_\infty)^{-1}\cdot \left( \|D_{\pi(\sp-1)} \bar \hol  D_{\pi(\sp-1)}\| + \| \bar\hol D_{\pi(\sp-1)} \odot \Sigma_{\sp-2} \| \right)\\
&\leq (1-\|\pi(\sp)\|_\infty)^{-1}\cdot \left(\|\pi(\sp)\|_\infty \|D_{\pi(\sp)} \hol  D_{\pi(\sp)}\| + \| D_{\pi(\sp)} \hol D_{\pi(\sp)} \odot \Sigma_{\sp-2} \| \right).
\end{align*}
Substituting the inequality above into \eqref{bound_D_pi_left} yields
\begin{align*}
 \| \hol \odot \Sigma_{\sp} \| &\leq (1-\|\pi(\sp)\|_\infty)^{-2}\cdot \left(\|D_{\pi(\sp)}  \hol  D_{\pi(\sp)}\| + \| D_{\pi(\sp)} \hol D_{\pi(\sp)} \odot \Sigma_{\sp-2} \| \right)
 \end{align*}
and since $\Sigma_{\sp-2}$ has $\pi(\sp-2)$ on its diagonal the result follows from Corollary~\ref{had} (with $A= \Sigma_{\sp-2}$) and Lemma~\ref{scalar_bounds}. \hfill {\bf (a)\checkmark} \\
\noindent \textbf{(b)} 
We want to show that $\E_\sp[\one_\Irv \one_\Irv^\tp] =\Sigma_{\sp} \preceq \pi \pi^\tp + 2 D_{\pi}$, meaning
\begin{equation*}
     \frac{\sum_{\Is: |\Is|=\sp} \one_\Is \one_\Is^\tp \PB(\Is)}{\sum_{\Is: |\Is|=\sp}\PB(\Is)} \preceq \frac{\sum_{(\Is,\Js): |\Is|=|\Js|=\sp} \one_\Is \one_\Js^\tp \PB(\Is)\PB(\Js)}{(\sum_{\Is: |\Is|=\sp}\PB(\Is))^2} +  \frac{2 \cdot\sum_{\Is: |\Is|=\sp} D_{\Is} \PB(\Is)}{\sum_{\Is: |\Is|=\sp}\PB(\Is)} .
\end{equation*}
Multiplying both sides by $(\sum_{\Is: |\Is|=\sp}\PB(\Is))^2$ we therefore have to show that
\begin{align*}
     \sum_{(\Is,\Js): |\Is|=|\Js|=\sp} \one_\Is \one_\Is^\tp \PB(\Is)\PB(\Js) &\preceq \sum_{(\Is,\Js): |\Is|=|\Js|=\sp} \one_\Is \one_\Js^\tp \PB(\Is)\PB(\Js)  + 2 \sum_{(\Is,\Js): |\Is|=|\Js|=\sp} D_{\Is} \PB(\Is)\PB(\Js).
\end{align*}
Now the crucial step is to partition these sums in a special way. For a pair $(\Is,\Js)$, by definition of the Poisson sampling model, we can rewrite $\PB(\Is)\PB(\Js)$ as
\begin{equation*}
    \PB(\Is)\PB(\Js) = \prod_{i \in \Is}p_i\prod_{j \notin \Is}(1-p_j)\prod_{i \in \Js}p_i\prod_{j \notin \Js}(1-p_j) = \prod_{i \in \Is \cap \Js}p_i^2\prod_{i \in \Is \triangle \Js}p_i(1-p_i)\prod_{j \notin \Is \cup \Js}(1-p_j)^2,
\end{equation*}
where $\Is \triangle \Js$ denotes the symmetric difference of $\Is,\Js$.
This implies that for two pairs $(\Is,\Js)$, $(\Is',\Js')$ we have
\begin{align*}
\quad \PB(\Is)\PB(\Js) = \PB(\Is')\PB(\Js')\quad 
\text{whenever}  \quad  \Is \cap \Js = \Is' \cap \Js' \quad \text{and} \quad  \Is \triangle \Js = \Is' \triangle \Js' .   
\end{align*}
This allows us to define natural partitions on the set of pairs $(\Is,\Js)$ such that the probability $\PB(\Is)\PB(\Js)$ is constant on each partition. Concretely, for any integer $\Tin \in \{1, \ldots , \sp\}$, together with a set $A\subseteq [K]$ with $|A| = \sp-\Tin$ and a set $B \subseteq [K]\setminus A $ with $|B| = 2\Tin$, we look at the collection of pairs $(\Is,\Js)$ with intersection $A$ and symmetric difference $B$, i.e.,
\begin{equation*}
    \mathcal{Q}_{A,B} := \left\{(\Is,\Js) : \Is,\Js \subseteq [K] , |\Is|=|\Js|=\sp , \Is \cap \Js =A, \Is\triangle \Js = B\right\}.
\end{equation*}
Note that each pair $(\Is,\Js)$ with $|\Is|=|\Js|=\sp$ can be \textit{uniquely} assigned to a collection $\mathcal{Q}_{A,B}$ and that $\PB(\Is)\PB(\Js)$ is constant for all $(\Is,\Js) \in \mathcal{Q}_{A,B}$. Using that the sum of positive semi-definite matrices is positive semi-definite, it therefore suffices to show that for all possible choices of $A,B$ we have
\begin{align*}
    \sum_{(\Is,\Js) \in \mathcal{Q}_{A,B}} \one_\Is \one_\Is^\tp \preceq \sum_{(\Is,\Js) \in \mathcal{Q}_{A,B}} \one_\Is \one_\Js^\tp  + 2 \sum_{(\Is,\Js) \in \mathcal{Q}_{A,B}} D_{\Is}. 
\end{align*}
For $A,B$ fixed we abbreviate $Q= \sum_{(\Is,\Js) \in \mathcal{Q}_{A,B}} \one_\Is \one_\Is^\tp$ and $\bar Q = \sum_{(\Is,\Js) \in \mathcal{Q}_{A,B}} \one_\Is \one_\Js^\tp$. Note that $\sum_{(\Is,\Js) \in \mathcal{Q}_{A,B}} D_{\Is}= Q \odot \Imat$, so the inequality above is equivalent to showing that
\begin{align}
0\preceq \bar Q - Q + 2 \cdot Q \odot \Imat.\label{poisson_d1}
\end{align}
For the entries of these matrices we have
\[
Q_{ij} = |\{ (\Is,\Js) \in \mathcal{Q}_{A,B} : i,j \in \Is \}| \quad\mbox{resp.}\quad
\bar Q_{ij} = |\{ (\Is,\Js) \in \mathcal{Q}_{A,B} : i\in \Is , j \in \Js \}|.
\]
In case $i,j\in A$ we obviously have $Q_{ij}=\bar Q_{ij} = |\mathcal{Q}_{A,B}|$, so
\[
Q \odot (\one_A \one_A^\tp) = \bar Q \odot (\one_A \one_A^\tp).
\] 
In particular, this means that \eqref{poisson_d1} holds trivially for $\Tin=0$, where $B=\emptyset$.\\
In case that $i \in A, j \in B$ we have
$Q_{ij}=\bar Q_{ij} = {2\Tin-1 \choose \Tin-1}=:d_\Tin$ and therefore
\[
Q \odot (\one_A \one_B^\tp + \one_B \one_A^\tp) = \bar Q \odot(\one_A \one_B^\tp + \one_B \one_A^\tp).
\]
It only remains to check what happens for $i,j\in B$. The case $\Tin = 0$ is already settled, thus we assume $\Tin\geq 1$. On the diagonal we get $\bar Q_{ii} =0$ while $Q_{ii} ={2\Tin-1 \choose \Tin-1} = d_\Tin$. We have $Q_{ij} = {2\Tin-2 \choose \Tin-2}=:q_\Tin$ and $\bar Q_{ij} = {2\Tin-2 \choose \Tin-1}=:\bar q_\Tin$. In summary
\begin{align}
\bar Q - Q + 2 \cdot Q \odot \Imat &= (\bar q_\Tin - q_\Tin) \cdot \one_B \one_B^\tp  - (\bar q_\Tin - q_\Tin + d_\Tin) \cdot D_B + 2 \cdot Q \odot \Imat. \label{poisson_d2}
\end{align}
Since $\bar q_\Tin \geq q_\Tin$ the matrix $(\bar q_\Tin - q_\Tin) \cdot \one_B \one_B^\tp$ is positive semi-definite. Finally, as $  {2\Tin-2 \choose \Tin-2} +  {2\Tin-2 \choose \Tin-1} =  {2\Tin-1 \choose \Tin-1}$ we have $\bar q_\Tin + q_\Tin = d_\Tin$. Thus for all $\Tin\geq 1$, we have
\begin{align*}
(\bar q_\Tin - q_\Tin + d_\Tin) \cdot D_B  &= 2\cdot  \bar q_\Tin \cdot D_B \preceq 2\cdot  d_\Tin \cdot  D_B  = 2 \cdot Q \odot D_B \preceq 2 \cdot Q \odot \Imat,
\end{align*}
showing that also the remaining terms in \eqref{poisson_d2} are positive semi-definite, which completes the proof of (b).
\hfill {\bf (b)\checkmark} \\
\noindent {\bf (c) }
 We will prove the statement by induction. Let $ \hat L$ be a set of size $\Tin\leq \sp-2$ and 
$\hat{\mathcal{L}} = \{\Is\subseteq [K]: \hat L \subseteq \Is\}$. We first show that for $k\notin \hat L$ and $L = \hat L \cup \{k\}$ we have
\begin{align}(1-\pi_k(\sp)) \cdot \E_{\sp} \big[ \one_{\Irv\setminus L} \one_{\Irv\setminus L }^\tp \cdot \indi_\mathcal{L}(\Irv) \big]  \preceq \pi_k(\sp) \cdot  \E_{\sp-1} \big[ \one_{\Irv\setminus \hat L} \one_{\Irv\setminus \hat L }^\tp \cdot \indi_{\hat{\mathcal{ L}}}(\Irv) \big]. \label{poisson_e1}
\end{align}
We again use that for any set $\Js$ not containing the index $k$ we have 
$$\frac{p_k}{1-p_k} \cdot \PB(\Js) =  \PB(\Js\cup\{k\}).$$
Thus expanding the expectation we get
\begin{align*}
\E_{\sp} \big[ \one_{\Irv\setminus L } \one_{\Irv\setminus L}^\tp \cdot \indi_{\mathcal{L}}(\Irv) \big]  
&=  \frac{\sum_{\Is :|\Is| = \sp, L \subseteq \Is} \PB(\Is)  (\one_{\Is\setminus L} \one_{\Is\setminus L}^\tp)}{\sum_{\Is : |\Is| = \sp } \PB(\Is)}  \cdot  \frac{\sum_{\Is: |\Is| = \sp , k \in \Is} \PB(\Is)}{\sum_{\Is : |\Is| = \sp , k \in \Is} \PB(\Is)}\nonumber \notag \\
&=  \frac{\sum_{\Is :|\Is| = \sp, L \subseteq \Is} \PB(\Is)  (\one_{\Is\setminus L} \one_{\Is\setminus L}^\tp)}{\sum_{\Is : |\Is| = \sp , k \in \Is} \PB(\Is)} 
\cdot  \frac{\sum_{\Is: |\Is| = \sp , k \in \Is} \PB(\Is)}{\sum_{\Is : |\Is| = \sp } \PB(\Is)}\nonumber \\
& =  \frac{\sum_{\Js : |\Js| = \sp-1, k\notin \Js, \hat L \subseteq \Js } \PB(\Js)  (\one_{\Js\setminus \hat L} \one_{\Js \setminus \hat L}^\tp)}{\sum_{\Js : |\Js| = \sp-1 , k \notin \Js} \PB(\Js)} \cdot \pi_k(\sp)  \nonumber \\
& \preceq   \frac{\sum_{\Js : |\Js| = \sp-1, \hat L \subseteq \Js } \PB(\Js)  (\one_{\Js\setminus \hat L} \one_{\Js \setminus \hat L}^\tp)}{\sum_{\Js : |\Js| = \sp-1 , k \notin \Js} \PB(\Js)} \cdot \frac{\sum_{\Is : |\Is| = \sp-1 } \PB(\Is)}{\sum_{\Is : |\Is| = \sp-1 } \PB(\Is)} \cdot \pi_k(\sp) \nonumber \\
&=   \frac{\sum_{\Js : |\Js| = \sp-1, \hat L \subseteq \Js } \PB(\Js)  (\one_{\Js\setminus \hat L} \one_{\Js \setminus \hat L}^\tp)}{\sum_{\Is : |\Is| = \sp-1 } \PB(\Is)} \cdot \frac{\sum_{\Is : |\Is| = \sp-1 } \PB(\Is)}{\sum_{\Js : |\Js| = \sp-1 , k \notin \Js} \PB(\Js)} \cdot \pi_k(\sp)  \nonumber \\
& =  \E_{\sp-1} \big[ \one_{\Irv \setminus \hat L} \one_{\Irv\setminus \hat L }^\tp \cdot \indi_{\hat{\mathcal{ L}}}(\Irv) \big] \cdot \frac{\sum_{\Is : |\Is| = \sp-1 } \PB(\Is)}{\sum_{\Js : |\Js| = \sp-1 , k \notin \Js} \PB(\Js)} \cdot \pi_k(\sp).
\end{align*}
Now all that remains to do in order to prove \eqref{poisson_e1} is to bound the fraction above. Writing out the expression in the denominator we get
\begin{align*}
\frac{\sum_{\Is : |\Is| = \sp-1 } \PB(\Is)}{\sum_{\Js : |\Js| = \sp-1 , k \notin \Js} \PB(\Js)}
&=  \frac{\sum_{\Is : |\Is| = \sp-1 } \PB(\Is)}{\sum_{\Is : |\Is| = \sp-1 } \PB(\Is) - \sum_{\Is : |\Is| = \sp-1 , k \in \Is} \PB(\Is)}  \\
&= \frac{1}{1- \P_{\sp-1}(k\in \Irv)} \leq \frac{1}{1- \P_{\sp}(k\in \Irv)}  = \frac{1}{1- \pi_k(\sp)}.
\end{align*}
By induction and using again that $\pi_k(\sp-1) \leq \pi_k(\sp)$, we finally get
$$\E_{\sp} \big[ \one_{\Irv\setminus L} \one_{\Irv\setminus L }^\tp \cdot \indi_{\mathcal{L}}(\Irv) \big] \prod_{\ell \in L} (1-\pi_{\ell}(\sp))  \preceq  \E_{\sp-|L|} [ \one_\Irv \one_\Irv^\tp] \cdot \prod_{\ell \in L} \pi_{\ell}(\sp), $$ 
which completes the proof of (c). \hfill \textbf{(c)\checkmark} \\
\phantom{bla for square position}
\end{proof}
Finally, it only remains to prove Corollary~\ref{cor_bound_HHt}.
\begin{proof}[of Corollary~\ref{cor_bound_HHt}]
We want to show that for a matrix $\hol$ with zero diagonal we have
\begin{align*}
\| \E[ D_{\Irv}\hol D_{\Irv} \hol^\tp D_{\Irv} ]\|  \leq& \frac{1 + \|\pi\|_\infty}{(1 - \|\pi\|_\infty)^3}  \|D_\pi H D_\pi H^\tp D_\pi \| +  \frac{1}{1 - \|\pi\|_\infty} \|D_\pi \odot (H D_\pi H^\tp) \| .
\end{align*}
Let $H_k$ be the $k$-th column of $H$ and $H_k^\tp = (H_k)^\tp$. We rewrite the expression, whose expectation we want to estimate, as 
\begin{align*}
D_{\Irv}\hol D_{\Irv} \hol^\tp D_{\Irv}& = (\hol  D_{\Irv} \hol^\tp) \odot (\one_{\Irv}\one_{\Irv}^\tp) = \big( \sum_{k \in \Irv} \hol_k \hol_k^\tp \big) \odot (\one_{\Irv}\one_{\Irv}^\tp)
  = \sum_{k \in \Irv} (\hol_k \hol_k^\tp) \odot (\one_{\Irv}\one_{\Irv}^\tp).
\end{align*}
Since the $k$-th entry of $\hol_k$ and therefore both the $k$-th row and $k$-th column of $\hol_k \hol_k^\tp$ are zero, we have $(\hol_k \hol_k^\tp) \odot (\one_{\Irv}\one_{\Irv}^\tp) = (\hol_k \hol_k^\tp) \odot (\one_{\Irv\setminus\{k\}} \one_{\Irv\setminus\{k\}}^\tp)$, yielding
\begin{align*}
D_{\Irv}\hol D_{\Irv} \hol^\tp D_{\Irv}& = \sum_{k \in \Irv} (\hol_k \hol_k^\tp) \odot (\one_{\Irv\setminus\{k\}} \one_{\Irv\setminus\{k\}}^\tp) \notag \\
 &= \sum_{k} (\one_{\Irv}(k) \cdot \hol_k \hol_k^\tp) \odot (\one_{\Irv\setminus\{k\}} \one_{\Irv\setminus\{k\}}^\tp)\notag \\
 &=\sum_{k}  (\hol_k \hol_k^\tp) \odot (\one_{\Irv\setminus\{k\}} \one_{\Irv\setminus\{k\}}^\tp \cdot \one_{\Irv}(k) ). 
\end{align*}
Note that $\one_{\Irv}(k) = \indi_{\mathcal{L}} (\Irv)$ for $L= \{k\}$ and $\mathcal{L} = \{ \Is \subseteq [K] : L \subseteq  \Is\}$. Thus using the Schur Product Theorem, which says that the Hadamard product of positive semi-definite (p.s.d) matrices is p.s.d., meaning for $A, P, \bar P$ p.s.d. with $P \preceq \bar P$ we have $A \odot P \preceq A \odot \bar P$, together with Theorem~\ref{matrix_bounds}\eqref{matrix_c} further leads to
\begin{align}
\E_{\sp} \left[ D_{\Irv}\hol D_{\Irv} \hol^\tp D_{\Irv} \right]&=
    \sum_{k} (\hol_k \hol_k^\tp) \odot \E_{\sp} \big[\one_{\Irv\setminus\{k\}} \one_{\Irv\setminus \{k\}} ^\tp \cdot \one_{\Irv}(k) \big] \notag\\
    &\preceq \sum_{k}  ( \hol_k   \hol_k^\tp ) \odot \big(\Sigma_{\sp-1} \cdot \tfrac{\pi_k}{1-\pi_k}\big) \notag \\
    &\preceq \underbrace{(1- \|\pi\|_\infty)^{-1}}_{=:c_\pi}  \big( \sum_{k}  \hol_k \cdot \pi_k \cdot \hol_k^\tp \big)  \odot \Sigma_{\sp-1} = c_\pi \cdot (\hol D_{\pi} \hol^\tp)\odot \Sigma_{\sp-1} .\notag 
\end{align}
This means that for the operator norm of the expression above we have the bound
\begin{align*}
\| \E_{\sp} [ D_{\Irv}\hol D_{\Irv} \hol^\tp D_{\Irv} ] \| &\leq c_\pi \cdot  \|  (\hol D_{\pi} \hol^\tp) \odot \Sigma_{\sp-1} ] \|,
\end{align*}
and all that remains to be done is to apply Theorem~\ref{matrix_bounds}\eqref{matrix_a} to $A = \hol D_{\pi} \hol^\tp$ in the bound above.
 So defining $C_{\pi(\sp)} :=   (1 + \|\pi(S)\|_\infty)(1 - \|\pi(S)\|_\infty)^{-2}$ and observing that due to Lemma~\ref{scalar_bounds} we have $\pi_k(\sp-1) \leq \pi_k(\sp) = \pi_k$ and thus $C_{\pi(\sp-1)} \leq C_{\pi(\sp)} =C_\pi$ we get 
 \begin{align*}
 \| \E_{\sp} [ D_{\Irv}\hol D_{\Irv} \hol^\tp D_{\Irv} ] \|  &\leq c_\pi \cdot  \| A \odot \Sigma_{\sp-1} ] \| \\
& \leq c_\pi \cdot \left( C_{\pi(\sp-1)} \cdot \|D_{\pi(\sp-1)} [A - A\odot \Imat\, ] D_{\pi(\sp-1)}\| +  \|A\odot D_{\pi(\sp-1)}\|\right) \\
& \leq c_\pi \cdot C_{\pi(\sp)} \cdot \|D_{\pi(\sp)} A D_{\pi(\sp)}\| +  c_\pi \cdot\| A \odot D_{\pi(\sp)}\|\\
& = c_\pi \cdot C_{\pi}\cdot \|D_{\pi} \hol D_{\pi} \hol^\tp D_{\pi}\| +  c_\pi \cdot \| D_\pi \odot (\hol D_{\pi} \hol^\tp)\|,
\end{align*}
where for the last inequality we have also used that $A=\hol D_{\pi}  \hol^\tp$ is positive semidefinite. 
\end{proof}

\appendix

\section{Motivation and derivation of Relation~\eqref{eq_newdico} } \label{sec_dico_update}
For the interested reader we give a short motivation how one constructs the glimpses in dictionary learning, which in first order approximation lead to Relation~\eqref{eq_newdico}. We then derive a similar relation using second order approximations and show how Corollary~\ref{cor_bound_HHt} comes into play.\\
The goal of dictionary learning is to identify the dictionary $\dico$ from a set of sparse signals, following a model $y=\dico_{\Irv} x_{\Irv}= \dico D_\Irv x$. The main idea behind most iterative dictionary learning algorithms is to take a guess $\pdico$ of the same size as $\dico$, and from each signal realisation $y =\dico_{\Is} x_{\Is}$ try to identify $\Is$ and subsequently $x_{\Is}$, using some sparse approximation scheme. Combinations of $y$ and $x_{\Is}$ are then aggregated in a way suitable to approximate $\dico$. For instance if ${\Irv}$ follows a rejective sampling model and for simplicity $x= \sigma$, a Rademacher sequence, we have 
$\E [yx^\tp] = \E [\dico D_\Irv x x^\tp]=  \E [\dico D_\Irv \Imat ] = \dico D_\pi$, so after renormalisation we can retrieve $\dico$ from $\E [yx^\tp]$ or rather its approximation by finite averages. Of course, if we only have $\pdico$ instead of $\dico$ we will not get $x$ but only an estimate $\hat x$ and the question is if using $\E[y\hat x^\tp]$ will lead to a better approximation of $\dico$ than $\pdico$.
Let us for simplicity assume that the sparse approximation routine always manages to retrieve the correct support realisation $\Is$, which is true with high enough probability that the resulting error remains small whenever $\pdico$ is close to $\dico$, cp. \cite{rusc26}.
Then the best approximation to $y$ using $\pdico_{\Is}$ is $\hat y  = \pdico_{\Is} \pdico_{\Is}^\dagger y$, meaning $\hat{x}_{\Is} = \pdico^\dagger_{\Is}y= \pdico^\dagger_{\Is} \dico_I x_{\Is}$, where $\pdico^\dagger_{\Is}$ is the Moore-Penrose pseudo-inverse. If $\pdico_{\Is}$ has full rank $\sp$, then we have
\begin{align*} 
\hat{x}_{\Is} =  \pdico^\dagger_{\Is} \dico_I x_{\Is} 
= \pdico^\dagger_{\Is} \pdico_{\Is} x_{\Is}  + \pdico^\dagger_{\Is} (\dico_{\Is} -\pdico_{\Is}) x_{\Is} 
= x_{\Is} + \pdico^\dagger_{\Is} Z_I x_{\Is} 
= x_{\Is} +  (\pdico_{\Is}^\tp \pdico_{\Is})^{-1} \pdico_{\Is}^\tp Z_I x_{\Is} .
\end{align*} 
If $G_{\Is,\Is} =  \pdico_{\Is}^\tp \pdico_{\Is} - \Imat $ has small operator norm, which for realisations of rejective sampling happens with high probability as long as $G = \pdico^\tp \pdico - \Imat$ has small entries and small weighted operator norm, \cite{rusc21}, then using a Neumann series we have 
$(\pdico_{\Is}^\tp \pdico_{\Is})^{-1} = \sum_{n=0}^\infty (- G_{\Is,\Is})^n$. This means that in first order we can approximate $ \pdico^\dagger_{\Is}$ by $ \pdico_{\Is}^\tp$ and in second order by $ \pdico_{\Is}^\tp -  G_{\Is,\Is}  \pdico_{\Is}^\tp $. So for most support realisations $\Is$ we have 
\begin{align*}
&\hat{x}_{\Is}\,\dot{=} \, x_{\Is} + \pdico_{\Is}^\tp Z_I x_{\Is} = \DI x + \DI  \pdico^\tp Z \DI x \quad \mbox{and}\\
&\hat{x}_{\Is} \,\ddot{=} \, x_{\Is} + \pdico_{\Is}^\tp Z_I x_{\Is} - G_{\Is,\Is}  \pdico_{\Is}^\tp Z_I x_{\Is} = \DI x + \DI  \pdico Z \DI x - \DI G \DI \pdico^\tp Z \DI x.
\end{align*}
Assuming again that $x= \sigma$ is a Rademacher sequence, we now arrive at Relation~\eqref{eq_newdico},
\begin{align*}
\E [y \hat x^\tp]\,\dot{=} \, \dico \, \E[  D_\Irv xx^\tp D_\Irv + D_\Irv xx^\tp D_\Irv Z^\tp \pdico D_\Irv] =  \dico [D_\pi + (Z^\tp \pdico)\odot \Sigma_\sp].
\end{align*}
Similarly, using the second order approximation and recalling that $Z^\tp \pdico = D_\eps + H$ we get
\begin{align*}
\E [y \hat x^\tp]\,\ddot{=} \,& \dico  \E[ D_\Irv +  D_\Irv Z^\tp \pdico D_\Irv - \DI Z^\tp \pdico D_\Irv G D_\Irv ] \\
=\,& \dico [D_\pi + (Z^\tp \pdico - D_\eps G)\odot \Sigma_\sp] -  \dico \E[ D_\Irv H D_\Irv G D_\Irv ] .
\end{align*}
Theorem~\ref{matrix_bounds} already allows us to directly bound the first two terms. To control the last term note that $\Is$ is a discrete random variable so we can interpret the expectation as a sum of matrices $\sum_\ell A_\ell B_\ell$. Let $A = (A_1,\ldots, A_L)$ and $B = (B_1^\tp, \ldots , B_L^\tp)^\tp$, then $\| \sum_\ell A_\ell B_\ell\|^2 = \| AB\|^2 \leq \|A\|^2 \cdot \|B\|^2 = \|AA^\tp\| \cdot \|B^\tp B\|$. Applying this to the expectation above where $A_\ell \sim \P_\sp(S=I)^{1/2} \DI H \DI$ and $B_\ell \sim \P_\sp(S=I)^{1/2} \DI G \DI$ leads to 
\begin{align*}
\|\E[ D_\Irv H D_\Irv G D_\Irv] \|^2 \leq \| \E[  D_\Irv H D_\Irv H^\tp D_\Irv ]\| \cdot \|\E[ D_\Irv G^\tp D_\Irv G D_\Irv ] \|,
\end{align*}
which we can finally bound using Corollary~\ref{cor_bound_HHt}. Note, that using tricks similar to the ones above together with Lemma~\ref{scalar_bounds}, Theorem~\ref{matrix_bounds} and Corollary~\ref{cor_bound_HHt} also allows to treat $\E [y \hat x^\tp]$ in full accuracy. This opens up the way to the convergence analysis of practically relevant dictionary learning algorithms, which usually involve additional components such as $\E[\hat x \hat x^\tp]$, see e.g \cite{rusc26}.

\section{Elementary proof of Lemma~\ref{scalar_bounds}} \label{sec_proof_scalar_bounds}
For the reader's convenience we here prove Lemma~\ref{scalar_bounds} using only the language and techniques introduced in the paper. 
\begin{proof}[of Lemma~\ref{scalar_bounds}]
We first want to show that for two indices $i\neq j$ we have
\begin{align*}
   \pi_{ij}(\sp) &=  \pi_i(\sp) \cdot \frac{\pi_j(\sp-1) -    \pi_{ij}(\sp-1)}{1 - \pi_i(\sp-1)}.
\end{align*}
Recall that for any set $\Js$ not containing the index $i$ we have 
$$\frac{p_i}{1-p_i} \cdot \PB(\Js) =  \PB(\Js\cup\{i\}).$$ Using this relation we get
\begin{align*}
   \pi_{ij}(\sp) &= \frac{\sum_{\Is:|\Is| = \sp} \indi_\Is(i) \indi_\Is(j) \cdot \PB(\Is)}{\sum_{\Is:|\Is| = \sp} \PB(\Is)} \cdot \frac{\sum_{\Is:|\Is| = \sp, i \in \Is} \PB(\Is)}{\sum_{\Is: |\Is| = \sp, i \in \Is} \PB(\Is)} \\
    & = \frac{\sum_{\Is:|\Is| = \sp, i \in \Is} \indi_\Is(j)\cdot \PB(\Is) } {\sum_{|\Is| = \sp, i \in \Is} \PB(\Is)} \cdot  \pi_i(\sp) \\
    & = \frac{p_i}{1-p_i} \cdot  \frac{1 -p_i}{p_i} \cdot \frac{\sum_{\Js:|\Js| = \sp-1, i \notin \Js} \indi_\Js(j) \cdot\PB(\Js) } {\sum_{\Js:|\Js| = \sp-1, i \notin \Js} \PB(\Js)}\cdot  \pi_i(\sp) \cdot \frac{\sum_{\Js:|\Js| = \sp-1} \PB(\Js)}{\sum_{\Js:|\Js| = \sp-1} \PB(\Js)}\\
    &= \pi_i(\sp) \cdot \frac{\sum_{\Js:|\Js| = \sp-1, i \notin \Js} \indi_\Js(j)\cdot \PB(\Js)} {\sum_{\Js:|\Js| = \sp-1} \PB(\Js)} \cdot \frac{\sum_{\Js:|\Js| = \sp-1} \PB(\Js)}{\sum_{\Js:|\Js| = \sp-1, i \notin \Js} \PB(\Js)}.
    \end{align*}
Further rewriting the fractions in the expression above yields 
    \begin{align*}
    \frac{\sum_{\Js:|\Js| = \sp-1, i \notin \Js} \indi_\Js(j) \cdot \PB(\Js)}{\sum_{\Js:|\Js| = \sp-1} \PB(\Js)}  
 &=  \underbrace{\frac{\sum_{\Js:|\Js| = \sp-1} \indi_\Js(j) \cdot \PB(\Is)} {\sum_{\Js:|\Js| = \sp-1} \PB(\Js)}}_{\pi_{j}(\sp-1)} - \underbrace{\frac{\sum_{\Js:|\Js| = \sp-1} \indi_\Js(i) \indi_\Js(j)  \cdot \PB(\Js) } {\sum_{\Js:|\Js| = \sp-1} \PB(\Js)} }_{   \pi_{ij}(\sp-1) } 
\end{align*}
as well as
\begin{align*}
    \frac{\sum_{\Js:|\Js| = \sp-1} \PB(\Js)}{\sum_{\Js:|\Js| = \sp-1, i \notin \Js} \PB(\Js)} &= \frac{\sum_{\Js:|\Js| = \sp-1} \PB(\Js)}{\sum_{\Js:|\Js| = \sp-1} \PB(\Js)- \sum_{\Js:|\Js| = \sp-1, i \in \Js} \PB(\Js)} = \frac{1}{1- \pi_{i}(\sp-1)},
\end{align*}
which completes the proof of the first equality.\\
\noindent Next we will show more generally that for arbitrary $L\subseteq [K]$, $\pi_L(\sp-1) \leq  \pi_L(\sp)$. We define $\mathcal L = \{ \Is\subseteq [K]: L \subseteq \Is\}$. Using this together with the definition of $\P_\sp$ we can rewrite $\pi_L(\sp-1) \leq  \pi_L(\sp)$ as
\begin{align*}
    \frac{ \sum_{\substack{\Js:|\Js| = \sp-1 }} \indi_\mathcal{L}(\Js) \cdot \PB(\Js) }{\sum_{\substack{\Js:|\Js| = \sp -1 }} \PB(\Js)} \leq \frac{ \sum_{\Is:|\Is| = \sp}\indi_\mathcal{L}(\Is) \cdot \PB(\Is)}{\sum_{\Is:|\Is| = \sp } \PB(\Is)},
\end{align*}
which is equivalent to
\begin{align*}
    \sum_{\substack{(\Is,\Js): |\Js| = \sp-1 , |\Is| = \sp }}\indi_\mathcal{L}(\Js)\cdot \PB(\Js)\PB(\Is) \leq \sum_{\substack{(\Is,\Js): |\Js| = \sp-1, |\Is| = \sp  }}\indi_\mathcal{L}(\Is)\cdot \PB(\Js)\PB(\Is).
\end{align*}
Next we use a similar partition as in the proof of Theorem~\ref{matrix_bounds}(b), that is for $\Tin \in \{1, \ldots , \sp\}$, together with a set $A\subseteq [K]$ with $|A| = \sp-\Tin$ and a set $B \subseteq [K]\setminus A $ with $|B| = 2\Tin-1$, we look at the collection of pairs $(\Is,\Js)$ with intersection $A$ and symmetric difference $B$, i.e.,
\begin{equation*}
    \mathcal{Q}_{A,B} := \left\{(\Is,\Js) : \Is,\Js \subseteq [K] , |\Is|=\sp, |\Js|=\sp-1 , \Is \cap \Js =A, \Is\triangle \Js = B\right\}.
\end{equation*}
Since each pair $(\Is,\Js)$ with $|\Is|=\sp, |\Js|= \sp-1$ can be uniquely assigned to a collection $\mathcal{Q}_{A,B}$ and $\PB(\Is)\PB(\Js)$ is constant for all $(\Is,\Js) \in \mathcal{Q}_{A,B}$, it is sufficient to show that
\begin{align*}
    \sum_{(\Is,\Js) \in \mathcal{Q}_{A,B}} \indi_\mathcal{L}(\Js) \leq \sum_{(\Is,\Js) \in \mathcal{Q}_{A,B}} \indi_\mathcal{L}(\Is)
\end{align*}
or equivalently that
\begin{align}
    |\{(\Is,\Js) \in \mathcal{Q}_{A,B} : L \subseteq \Js \}|  \leq |\{(\Is,\Js) \in \mathcal{Q}_{A,B} : L \subseteq \Is \}. \label{pi_\sp_comp_a}
\end{align}
If $L$ is not contained in $A\cup B$, there is no valid pair $(\Is,\Js) \in \mathcal{Q}_{A,B}$ and the inequality trivially holds. 
If $L \subseteq A\cup B$ we abbreviate $L_A = L \cap A$ and $L_B = L\cap B$. Since $L_A \subseteq A$, all pairs
in $(\Is,\Js) \in \mathcal{Q}_{A,B} $ automatically satisfy $L_A \subseteq \Is$ and $L_A \subseteq \Js$ so 
we can rewrite~\eqref{pi_\sp_comp_a} as
\begin{align}
 |\{(\Is,\Js) \in \mathcal{Q}_{A,B}:  L_B \subseteq \Js \}|  \leq |\{(\Is,\Js) \in \mathcal{Q}_{A,B} : L_B \subseteq \Is\} |.\label{pi_\sp_comp_b}
\end{align}
Since we need to have $A \cup L_B \subseteq \Js$, in case $|A \cup L_B| = |A| + |L_B| \geq \sp$ the left hand side in \eqref{pi_\sp_comp_b} is zero and the inequality holds. Finally, if $k=|L_B|  < \sp - |A| = \Tin$, we can still choose
$\Tin-k-1$ out of the $2\Tin-k-1$ resp. $\Tin-k $ out of the $2\Tin-k  -1 $ remaining elements in $B$ to fill $\Is$ resp. $\Js$ and create a valid pair. Since 
\begin{align*}
 \binom{2\Tin-k-1}{\Tin- k-1} \leq  \binom{2\Tin-k -1}{\Tin-k  } ,
\end{align*}
the inequality in \eqref{pi_\sp_comp_b} is satisfied, which completes the proof.
\end{proof}
Finally, let us quickly exploit the partitioning trick to get an elementary proof for 
another bound, that comes in handy when working with rejective sampling inclusion probabilities,
\begin{align}
\pi_{L \cup M}(\sp) \leq \pi_L(\sp) \cdot \pi_M(\sp) \quad \mbox{if} \quad M\cap L = \emptyset.
\end{align}
Again this bound can be proved using 3 lines and some heavy machinery from \cite{neg_dep_09}, i.e., since $\mu = \PB$ is strongly Raleigh, so is $\mu_n = \P_n$, which then  implies that $\P_n$ is negatively associated, so for increasing functions $f = \sum_{\ell \in L} \one_\Irv(\ell)$ and $g =\sum_{m\in M} \one_\Irv(m)$ we have 
$\pi_{L \cup M}(\sp) = \int fg \mathrm{d}\mu_n \leq  \int f \mathrm{d}\mu_n  \int g \mathrm{d}\mu_n =\pi_L(\sp) \cdot \pi_M(\sp)$. \\
Alternatively, we can use somewhat more lines but only elementary combinatorics.\\
We define $\mathcal L = \{ \Is\subseteq [K]: L \subseteq \Is\}$ and $\mathcal M = \{ \Is\subseteq [K]: M \subseteq \Is\}$. Using this together with the definition of $\P_\sp$ we can rewrite $\pi_{L \cup M}(\sp) \leq \pi_L(\sp) \cdot \pi_M(\sp)$ as
\begin{align*}
    \frac{ \sum_{\substack{\Is:|\Is| = \sp }} \indi_\mathcal{M}(\Is)\cdot \indi_\mathcal{L}(\Is) \cdot \PB(\Is) }{\sum_{\substack{\Is:|\Is| = \sp }} \PB(\Is)} \leq  \frac{ \sum_{\substack{\Js:|\Js| = \sp }}\indi_\mathcal{M}(\Js)\cdot \PB(\Js)}{\sum_{\substack{\Js:|\Js| = \sp }} \PB(\Js)} \cdot \frac{ \sum_{\Is:|\Is| = \sp}\indi_\mathcal{L}(\Is) \cdot \PB(\Is)}{\sum_{\Is:|\Is| = \sp } \PB(\Is)},
\end{align*}
which is equivalent to
\begin{align*}
    \sum_{(\Is,\Js):|\Is| = |\Js| = \sp} \indi_\mathcal{M}(\Is) \cdot \indi_\mathcal{L}(\Is)\cdot \PB(\Is)\PB(\Js) 
    \leq \sum_{(\Is,\Js):|\Is| = |\Js| = \sp} \indi_\mathcal{M}(\Js) \cdot  \indi_\mathcal{L}(\Is) \cdot \PB(\Is)\PB(\Js).
\end{align*}
We now use the same decomposition as in the proof of Theorem~\ref{matrix_bounds}(b), meaning for $\Tin \in \{0, \ldots , \sp\}$, $A\subseteq [K]$ with $|A| = \sp-\Tin$ and $B \subseteq [K]\setminus A $ with $|B| = 2\Tin$, we define
\begin{equation*}
    \mathcal{Q}_{A,B} := \left\{(\Is,\Js) : \Is,\Js \subseteq [K] , |\Is|= |\Js|=\sp, \Is\cap \Js =A, \Is\triangle \Js = B\right\}.
\end{equation*}
Using the usual argument we now only need to show that 
\begin{align*}
    \sum_{(\Is,\Js) \in \mathcal{Q}_{A,B}} \indi_\mathcal{M}(\Is) \cdot \indi_\mathcal{L}(\Is) \leq \sum_{(\Is,\Js) \in \mathcal{Q}_{A,B}} \indi_\mathcal{M}(\Js) \cdot  \indi_\mathcal{L}(\Is)
\end{align*}
or equivalently that 
\begin{align*}
    |\{(\Is,\Js) \in \mathcal{Q}_{A,B} :M \subseteq \Is, L \subseteq \Is \}|  \leq |\{(\Is,\Js) \in \mathcal{Q}_{A,B} : M \subseteq \Js, L \subseteq \Is\}|. 
\end{align*}
If $L \cup M$ is not contained in $A\cup B$ there is no valid pair $(\Is,\Js) \in \mathcal{Q}_{A,B}$ and the inequality trivially holds. 
If $(L \cup M) \subseteq A\cup B$ we abbreviate $L_A = L \cap A$, $L_B = L\cap B$, $M_A = M \cap A$ and $M_B = M\cap B$. Since $(L_A \cup M_A) \subseteq A$, all pairs
in $(\Is,\Js) \in \mathcal{Q}_{A,B} $ automatically satisfy $(L_A \cup M_A) \subseteq \Is$, $M_A \subseteq \Js$ and $L_A \subseteq \Is $ so 
we can rewrite the inequality we want to show as
\begin{align}
 |\{(\Is,\Js) \in \mathcal{Q}_{A,B} :M_B \subseteq \Is, L_B \subseteq \Is \}|  \leq |\{(\Is,\Js) \in \mathcal{Q}_{A,B} : M_B \subseteq \Js, L_B \subseteq \Is\} | . \label{poisson_c1}
\end{align}
Since we need to have $(A \cup L_B \cup M_B) \subseteq \Is$, in case $|A \cup L_B \cup M_B| = |A| + |L_B| + |M_B| > \sp$ the left hand side in \eqref{poisson_c1} is zero and the inequality trivially holds. Finally, if $k=|L_B| + |M_B|  \leq \sp - |A| = \Tin$, we can still choose
$\Tin-k$ out of the $2\Tin-k$ resp. $\Tin- |L_B|  $ out of the $2\Tin-k$ remaining elements in $B$ to fill $\Is$ resp. $\Js$ and create a valid pair. Since 
$|L_B|\leq k$, we have \begin{align*}
 \binom{2\Tin-k}{\Tin- k} \leq  \binom{2\Tin-k }{\Tin- |L_B | }, 
\end{align*}
meaning the inequality in \eqref{poisson_c1} is again satisfied, which completes the proof.

\acks{This work was supported by the Austrian Science Fund (FWF) under Grant no.~Y760.
}

\end{document}